\documentclass{amsart}
\usepackage{graphicx}

\theoremstyle{definition}

\theoremstyle{remark}

\numberwithin{equation}{section}

\begin{document}

\title{Improving Riemann prime counting}

\author{Michel Planat}
\address{Institut FEMTO-ST, CNRS, 15 B Avenue des Montboucons, F-25044 Besan\c con, France.}
\email{michel.planat@femto-st.fr}

\author{Patrick Sol\'e}
\address{Telecom ParisTech, 46 rue Barrault, 75634 Paris Cedex 13, France.}
\curraddr{Mathematics Department, King Abdulaziz University, Jeddah, Saudi Arabia.}
\email{sole@telecom-paristech.fr}


\date{October 1\normalfont, 2014}



\begin{abstract}

Prime number theorem asserts that (at large $x$) the prime counting function $\pi(x)$ is approximately the logarithmic integral $\mbox{li}(x)$. In the intermediate range, Riemann prime counting function $\mbox{Ri}^{(N)}(x)=\sum_{n=1}^N \frac{\mu(n)}{n}\mbox{Li}(x^{1/n})$ deviates from $\pi(x)$ by the asymptotically vanishing sum $\sum_{\rho}\mbox{Ri}(x^\rho)$ depending on the critical zeros $\rho$ of the Riemann zeta function $\zeta(s)$. We find a fit $\pi(x)\approx \mbox{Ri}^{(3)}[\psi(x)]$ [with three to four new exact digits compared to $\mbox{li}(x)$] by making use of the Von Mangoldt explicit formula for the Chebyshev function $\psi(x)$. Another equivalent fit makes use of the Gram formula with the variable  $\psi(x)$. Doing so, we evaluate $\pi(x)$ in the range $x=10^i$, $i=[1\cdots 50]$ with the help of the first $2\times 10^6$ Riemann zeros $\rho$. A few remarks related to Riemann hypothesis (RH) are given in this context.

\end{abstract}

\maketitle


\section{Introduction}

In his celebrated 1859 note about the prime counting function $\pi(x)$, Riemann concludes

 {\it The thickening and thinning of primes which is represented by the periodic term in the formula has also been observed in the counts of primes, without, however, any possibility of establishing a law for it having been noticed. It would be interesting in a future count to examine the influence of individual periodic terms in the formula for the density of primes \cite[p. 305]{Edwards74}}.

It is known that the periodic terms $\mbox{Ri}(x^\rho)$ at the critical zeros $\rho$ of the Riemann zeta function $\zeta(s)$ are responsible for the inaccuracy of Riemann prime counting function $\mbox{Ri}^{(N)}(x)=\sum_{n=1}^N \frac{\mu(n)}{n}\mbox{Li}(x^{1/n})$ (where $\mu(n)$ is the Moebius function) in the approximation of $\pi(x)$, the number of primes up to $x$ \cite[eq. (4)]{Borwein} 
$$\pi(x) \sim \pi_0(x)= \mbox{Ri}^{(N)}(x)-\sum_{\rho}\mbox{Ri}(x^\rho).$$

But even the account of the first $200$ periodic terms performed in the approximation of $\pi(10^{20})$ is not able to improve the approximation $\mbox{li}(10^{20})$ \footnote{In the following when we write $\mbox{li}$, we mean its integer part $\left\lfloor \mbox{li}\right\rfloor$.} by more that one digit \cite[p. 249]{Borwein}
$$\pi(10^{20})=u~560918840,$$
$$\mbox{li}(10^{20})=u~783663484,$$
$$\pi_0(10^{20})\sim u~591885820,$$
where $u$ means $2220819602$. In contrast, in \cite{PlanatSole}, we introduced a different approach where the account of the periodic terms is performed in a  global way as $\pi(x)\sim \pi_1(x)=\mbox{Ri}^{(3)}[\psi(x)]$ so that one obtains
$$\pi(10^{20})\sim \pi_1(10^{20})=u~561025931,$$ an almost three digit improvement over the approximation by the logarithmic integral.

Similarly, we are able to improve Gram formula \cite[eq. (70)]{Borwein}
$$\mbox{Ri}(x)=1+\sum_{n=1}^\infty \frac{(\log~ x)^n}{n \zeta(n+1) n!}$$ from the approximation $\pi_2(x)=1+\mbox{Ri}[\psi(x)]$ so that 
$$\pi(10^{20})\sim \pi_2(10^{20})=u~561025719.$$

This improvement is observed at all values $x=10^i$, $i \le 25$, where $\pi(x)$ is exactly known (see Sec. 2). 

One is fortunate to have now at our disposal two distinct (but closely related and surprisingly almost equal) formulas $\pi_1(x)$ and $\pi_2(x)$ to approximate the prime counting function in the unknown range of values such as $10^i$, $i>25$. We can use either formula to explore the unknown range of values such as $x:=10^i$, $i>25$. Specifically 

$$\mbox{li}(10^{50})=v~ 862622818995697067491328 ,$$
$$\pi(10^{50})\sim \pi_1(10^{50})=v~780454103362367083511808,$$
$$\pi(10^{50})\sim \pi_2(10^{50})=v~780552101690258665504768,$$ 
where $v$ means $876268031750784168878176$.

In this short paper, we pursue the calculations of $\pi(x)\sim \pi_1(x)$ started in \cite{PlanatSole} and compare it to the Gram formula based approximation $\pi_2(x)$ by having recourse to the explicit Von Mangoldt formula for the Chebyshev function $\psi(x)$ and taking into account the first $2 \times 10^6$ critical zeros of the Riemann zeta function \cite{Odlyzko}. These calculations are well in the spirit of the aforementioned quote of Riemann since the approximate value of $\pi(x)$ explicitly depends on the zeros $\rho$ through $\psi(x)$. The function $\mbox{li}[\psi(x)]$ also relates to Riemann hypothesis through the modified Robin's criterion \cite[Corollary 2]{PlanatSole},\cite{Robin84} so that we are close to this important topic as well. 

\section {Approximation of $\pi(x)$ at powers $x=10^i$, $i=[1\cdots 25]$}

It is claimed in the introduction that $\pi_1(x)$ and $\pi_2(x)$ provide a much better fit (two to four digits better) than $\mbox{li(x)}$. It also provides a much better fit than the Gram formula for $\mbox{Ri}(x)$ as shown below.  We used the first $7 \times 10^5$ zeros $\rho$ calculated in \cite{Odlyzko} to evaluate $\psi(x)$ and the corresponding approximations $\pi_1(x)$ and $\pi_2(x)$ from the explicit Von Mangoldt formula
$$\psi_0(x)=x-\sum_{\rho} \frac{x^{\rho}}{\rho}-\frac{\zeta'(0)}{\zeta(0)}-\frac{1}{2}\log (1-x^{-2}), ~\mbox{for}~x>1,$$
$\psi(x)=\psi_0(x)$ when $x$ is not a prime power and $\psi_(x)=\psi_0(x)+\frac{1}{2}\Lambda(x)$, $\Lambda(x)$ the Von Mangoldt function, otherwise \cite[p. 104]{Davenport80}. \footnote{Incidently, we mention and numerically check a remarkable formula, known to Riemann 
\cite[p. 67]{Edwards74}, about the properties of critical zeros 
$$\sum_{\rho} \frac{1}{\rho(1-\rho)}=2+\gamma- \log(4 \pi) =C.$$
By using the first $2 \times 10^6$ zeros \cite{Odlyzko}, one gets  $\sum_{\rho} \frac{1}{\rho(1-\rho)}\sim 0.999920.$.

Of course, RH is equivalent to $\sum_{\rho} \frac{1}{|\rho|^2}=C$, as observed in \cite{Burnol98}.}

The difference $\mbox{Ri}(x)-\pi(x)$ at powers of $10^i$, $i=[1\cdots 25]$ is the sequence A215663 

\begin{eqnarray}
&\{0,0,0,-3,-5,29,88,96,-79,-1828,   -2319,-1476,-5774,-19201, 73217, \nonumber\\
&327052,-598255,-3501366,23884333,-4891825,       -86432205,-127132665,\nonumber \\&1033299853,-1658989720,-1834784715\}.\nonumber \\ 
\nonumber
\end{eqnarray}

The sequence $\pi_1(x)-\pi(x)$  reads
\begin{eqnarray}
&\{0, 0, 0, 0, 1, 1, -1, 4, 4, -2,   24, -16, 67, -273, -2886, \nonumber \\
 & -5203, 24767, 39982, 11722, 107091, -339757, -3972640,,\nonumber \\
& 8296270, 75888611, -42528602\}.\nonumber \\
\nonumber
\end{eqnarray}
As for the sequence $\pi_2(x)-\pi(x)$, one gets
\begin{eqnarray}
&\{0, 0, 0, 0, 0, 0, -2, 2, 1, -6,   18, -24, 58, -280, -2852, \nonumber \\
&-5390, 24170, 39667, 9990, 106880, -372719, -3896886,\nonumber \\
&8380617, 75606965, -10884280\}.\nonumber \\
\nonumber
\end{eqnarray}

Clearly the shift from $\pi(x)$ of the approximation $\pi_1(x)$ and $\pi_2(x)$ is lower, sometimes with two orders of magnitude, than the shift from $\mbox{Ri}(x)$ with respect to $\pi(x)$, as expected.

It turns out that $\pi_2(x)=\pi(x)$ at all powers of $10$ up to $10^6$. But there exists plenty values of $x<10^6$ where $\pi_2(x) \ne \pi(x)$. This occurs over thick intervals starting at squares of primes $p^2$ with $p \in \{13,19,23,31,47,53,61,71,73,79,83,89\cdots\}$. The first $7$ primes in this sequence are associated with prime gaps (OEIS sequence A134266) but the rest is not recognized and quite random.

\section {Approximation of $\pi(x)$ at powers $x=10^i$, $i=[26\cdots 50]$}

The values of $\pi(10^i)$ at $i>25$ are not known. An approach to approximate $10^{26}$ is given in \cite{Pletser13}. It is already interesting to guarantee which digits of $\mbox{li}(x)$ are exact and which ones can be added. Our approximation based on $\pi_1(x) \sim \pi_2(x)$ defined in the introduction allows to solve this challenge as soon as a sufficient number of critical zeros of $\zeta(s)$ are taken into account. As before, we used the first $2 \times 10^6$ zeros calculated in \cite{Odlyzko} to evaluate $\psi(x)$ and the corresponding approximations of $\pi(x)$. The results are listed below. In the array, the prefix that belongs to  $\mbox{li}(x)$ is denoted $w$ and three to four exact digits that do not belong to $\mbox{li(x)}$ can be guaranteed; for further comments compare the evaluation of $\pi(10^{50})$ shown in the introduction.

\newpage
\small
\begin{eqnarray}
& \mbox{li}(10^{26})=1699246750872~ 593033005722, ~~\pi(10^{26}) \sim w~4370***\nonumber \\
& \mbox{li}(10^{27})=1635246042684~ 2189113085404,   ~~\pi(10^{27}) \sim w~1681***\nonumber\\
& \mbox{li}(10^{28})=15758926927597~ 4838158399970,   ~~\pi(10^{28}) \sim w~3411***\nonumber\\
& \mbox{li}(10^{29})=152069810971427~ 6717287880526,   ~~\pi(10^{29}) \sim w~2164***\nonumber\\
& \mbox{li}(10^{30})=146923988977204~ 47639079087666,   ~~\pi(10^{30}) \sim w~2886***\nonumber\\
& \mbox{li}(10^{31})=142115097348080~ 932014151407888,   ~~\pi(10^{31}) \sim w~8927***\nonumber\\
& \mbox{li}(10^{32})=137611086699376~ 6004917522323376,   ~~\pi(10^{32}) \sim w~5914***\nonumber\\
& \mbox{li}(10^{33})=13338384833104449~ 976987996078656,   ~~\pi(10^{33}) \sim w~5799***\nonumber\\
& \mbox{li}(10^{34})=12940862650515894~ 2694528690220032,   ~~\pi(10^{34}) \sim w~0734***\nonumber\\
& \mbox{li}(10^{35})=125663532881831647~ 7984258713700352,   ~~\pi(10^{35}) \sim w~3088***\nonumber\\
& \mbox{li}(10^{36})=122129142976193654 ~~54914525219717120,   ~~\pi(10^{36}) \sim w~4747***\nonumber\\
& \mbox{li}(10^{37})=1187881589121682517~ 51964753226629120,   ~~\pi(10^{37}) \sim w~2948***\nonumber\\
& \mbox{li}(10^{38})=1156251261026516898~ 117818362099662848,   ~~\pi(10^{38}) \sim w~0207***\nonumber\\
& \mbox{li}(10^{39})=1126261940555920314~ 6275209660101296128,   ~~\pi(10^{39}) \sim w~5914***\nonumber\\
& \mbox{li}(10^{40})=10977891348982830282~ 8088665394982682624,   ~~\pi(10^{40}) \sim w~7136***\nonumber\\
& \mbox{li}(10^{41})=107072063488003546554~ 7951820678396641280,   ~~\pi(10^{41}) \sim w~5104***\nonumber\\
& \mbox{li}(10^{42})=104495503622645875354~ 81840562999989895168,   ~~\pi(10^{42}) \sim w~7283***\nonumber\\
& \mbox{li}(10^{43})=1020400469443659108805~ 59853621001121169408,   ~~\pi(10^{43}) \sim w~2642***\nonumber\\
& \mbox{li}(10^{44})=996973504768769817629~ 283320476467636207616,   ~~\pi(10^{44}) \sim w~1914***\nonumber\\
& \mbox{li}(10^{45})=9745982046649286035485~ 484938845939166085120,   ~~\pi(10^{45}) \sim w~2597***\nonumber\\
& \mbox{li}(10^{46})=9532053011747645833855~ 1157512823190583246848,   ~~\pi(10^{46}) \sim w~0125***\nonumber\\
& \mbox{li}(10^{47})= 93273147934738153021141 ~7209650549631557304320,   ~~\pi(10^{47}) \sim w~3599***\nonumber\\
& \mbox{li}(10^{48})= 91311875111614162331019~~ 53278507015309237420032,   ~~\pi(10^{48}) \sim w~4673***\nonumber\\
& \mbox{li}(10^{49})=894313906580259138316220~ 54467196576109747503104 ,   ~~\pi(10^{49}) \sim w~3207***\nonumber\\
& \mbox{li}(10^{50})=876268031750784168878176~ 862622818995697067491328 ,   ~~\pi(10^{50}) \sim w~7805***\nonumber\\
 \nonumber\\   
\nonumber
\end {eqnarray}
\normalsize

\section{Approximation of $\pi(x)$ where $\psi(x)$ is exactly known}

To check the validity of prime counting functions $\pi_1(x)$ and $\pi_2(x)$ it is good to compute them at values where the Chebyshev function $\psi(x)$ is exactly known, that is, irrespectively of the knowledge of the critical zeros $\rho$. Exact values of $\psi(x)$ at selected high values of $x$, with $10^6 \le x \le 10^{15}$ are given in \cite{Deleglise}. 

In this subsection we restrict to the calculation of the Gram formula based approximation $\pi_2(x)$. But similar observations hold for $\pi_1(x)$.

In the following two lists, $\pi'_2(x)$ is calculated from the (almost) exact values of $\psi(x)$ and $\pi_2(x)$ is calculated from the explicit formula with $2 \times 10^6$ critical zeros. The shifts $\pi'_2(x)-\pi(x)$  and $\pi_2(x)-\pi(x)$ are given at values of $x$ found in \cite[Table 1]{Deleglise}, that is, $x \in \{k. 10^j\}$, $k=[1\cdots 9]$, $j=[6\cdots 14]$ and at $x=10^{15}$. To facilitate the reading of the lists, we put a semi-column at $x$ values preceding every power of $10$.

\small
\begin{eqnarray}
& \pi'_2(x)-\pi(x)= \nonumber\\
&\{;0, -1, 1, 0, -1, 2, -1, 2, 0; -2, 0, 0, 1, 2, 1, 0, 0, 2;
2, -2, -1, -1, 0, 4, 1, -1, -2;~~ \nonumber\\
& -1, 3, -2, -3, -3, 0, 2, -2, 4; 2,8,-3,-2,8,-2,-5,6,-12;-10,11,8,-4,6,-14,12,16,-9;~~     \nonumber \\
& -23,16,-8,5,13,-21,8,-3,-17; -27,-24,5,76,15,66,28,-46,81; \nonumber \\
&-9,-132,-46,120,-65,302,-214,-11,197;168\}   \nonumber \\
\nonumber
\end{eqnarray}
\normalsize

\small
\begin{eqnarray}
& \pi_2(x)-\pi(x)= \nonumber\\
&\{; 0, -1, 0, 0, -1, 2, -1, 2, 0;-2, 0, 1, 2, 3, 1, -1, 0, 1;2, -1, -2, -2, 6, 7, -2, -2, 1;  \nonumber\\
&1, 7, -4, -7, 4, 10, 4, 19, -13; -6,17,-37,-15,-2,30,11,48,22;18,44,11,25,9,12,-59,36,-8;     \nonumber \\
&-24,-144,-34,-292,77,252,-81,-410,5;58,61,-6,58,258,-894,719,556,-401;
   \nonumber \\
& -280,94,-842,-1028,178,1425,597,247,-1617;-2852 \}   \nonumber \\
\nonumber
\end{eqnarray}
\normalsize

It is clear that, while the shifts are almost equal (and very small) at the beginning of the lists, they are higher at the end of the lists, and they differ substantially (about one order of magnitude) in the two. 

These calculations reinforce our confidence in the efficiency of the prime counting function $\pi_2(x)$ in that the remaining inaccuracy of $\pi_2(x)$ partially arises from the possible inacurracy of the calculation of $\psi(x)$ by the explicit formula. Unfortunately, at high values of $x$, the time for computing $\psi(x)$ becomes as prohibitive as the time for computing $\pi(x)$ \cite{Deleglise}, and this is why the approximation of $\psi(x)$ based on the explicit formula remains extremely useful.

Looking at the relative error $\epsilon:=[\pi_2(x)-\pi(x)]/\pi(x)$ compared to $\eta:=[\mbox{li}(x)-\pi(x)]/\pi(x)$, over the range of the above explored values $10^6<x<10^{25}$, one gets $0 \le |\epsilon/\eta|<6.5 \times 10^{-3}$ but the average ratio $|\epsilon/\eta|$ is about $1.2 \times 10^{-3}$. This represents an improvement of about three orders of magnitude of the prime counting function $\pi_2(x)$ compared to $\mbox{li}(x)$. Depending on the selected value of $x$, three to four new exact digits are obtained from $\pi_2(x)$ compared to  $\mbox{li}(x)$ as shown in the previous sections.

\section{Hints about the function $\mbox{li}(\psi(x))$ and RH}

Littlewood established that the function $\pi(x)-\mbox{li}(x)$ changes sign infinitely often. But it is known not to occur before $x=x_0\approx e^{727}$, a so-called Skewes' number \cite{Skewes}. 

Asymptotically, one has $x \sim\theta(x) \sim \psi(x)$, where $\theta(x)=\sum_{p\le x} \log p$ is the first Chebyshev function. But Robin proved the statement \cite{Robin84}
$$\epsilon_{\theta(x)}=\mbox{li}[\theta(x)]-\pi(x)>0 ~~\mbox{is}~~\mbox{equivalent}~~\mbox{to}~~\mbox{RH}.$$ 
As a corollary, the statement $\epsilon_{\psi(x)}=\mbox{li}(\psi(x))-\pi(x)>0$ is also equivalent to RH as was already observed in \cite{PlanatSole}. \footnote{In a related work, we introduce a similar statement \cite[eq. (2.1)]{ChebBias11} as a Chebyshev's type bias whose positivity is equivalent to GRH for the corresponding modulus.}

 As we arrived at the  excellent counting functions $\pi_1(x)$ and $\pi_2(x)$, themselves functions of $\psi(x)$ and thus explicitly related to the zeros of $\rho$ of the Riemann zeta function, it is quite satisfactory to be back to the spirit of Riemann's program of counting the prime numbers.

{\it One would of course like to have a rigorous proof of this, but I have put aside the search for such a proof after some fleeting vain attemps  because it is not necessary for the immediate objective of my investigation} \cite[p. 301]{Edwards74}.

To conclude, Riemann prime counting function $\mbox{Ri}(x)$ can be much improved by replacing the variable $x$ by the Chebyshev function $\psi(x)$, but it is challenging to understand the origin of this seemingly \lq\lq explicit formula" for $\pi(x)$.

\bibliographystyle{amsplain}


\end{document}